\documentclass[letter,11pt]{article}
\usepackage{amsfonts}
\usepackage{amssymb}
\usepackage{amsmath}
\usepackage{cite}
\newtheorem{theo}{Theorem}
\newtheorem{prop}[theo]{Proposition}
\newtheorem{rema}[theo]{Remark}
\newtheorem{lema}[theo]{Lemma}

\newtheorem{exam}[theo]{Example}

\newcommand{\R}{\mathbb{R}}
\newcommand{\D}{\mathbb{D}}
\newcommand{\T}{\mathbb{T}^1}

\newcommand{\C}{\mathbb{C}}
\newcommand{\N}{\mathbb{N}}
\newcommand{\eps}{\varepsilon}
\newcommand{\Z}{\mathbb{Z}}

\title{Naishul's theorem for fibred holomorphic maps}
\author{Mario PONCE\footnote{Partially supported by FONDECYT 11090003}\\
PUC-Chile\\
mponcea@mat.puc.cl}
\begin{document}
\maketitle
\begin{abstract}
We show that the fibred  rotation number associated to an indifferent invariant curve for a fibred holomorphic map is a topological invariant.
\end{abstract}

\section{Introduction}
The rotation number was introduced by Poincar\'e in order to compare a general circle transformation with the simplest non trivial dynamical model, the rigid rotations. Thus, any orientation preserving (positive) circle homeomorphism has a  rigid rotation  as a model for its dynamics.  Poincar\'e shows that this number characterizes the order of the orbits along the circle, controlling the topology or the {\it shape} of any orbit. By that, and a simple computation, the rotation number results into a topological conjugacy invariant. Moreover, under some arithmetical and smoothness conditions one can say that the rotation number is a characterization  of the full conjugacy class: these conditions  imply that the map is conjugated to the corresponding rigid rotation (Denjoy, Arnold, Herman, Yoccoz, see for instance \cite{KAHA95} and references therein). \\

A rotation number can be associated to a differentiable surface local homeomorphism having an indifferent fixed point, the so called {\it non linear rotations}. The derivative at the fixed point is a pure rotation and one figure that the rotation number have some control on the dynamics of points which are close to the fixed point. By assuming some hypothesis one can give results in this direction.    For example, in the holomorphic case, and under the Brjuno arithmetical condition, the dynamics is actually conjugated to the pure rotation map.  Even in the absence of this nice behavior, the rotation number determines the shape  of orbits in the following sense.
\begin{theo}[Naishul \cite{NAIS82}]
Let $f$ and $g$ be two orientation preserving non linear rotations, which are topologically conjugated by a conjugacy which preserves orientation and the fixed point. If $f$ is holomorphic (or area preserving), then the rotation numbers of $f$ and $g$ are equals $\quad_{\blacksquare}$
\end{theo}
In a subsequent work, Gambaudo and P\'ecou \cite{GAPE95} realized that the smoothness condition  (holomorphic or area preserving) is not the intrinsic ingredient which turns true this result. They define the {\it linking property} for nearby orbits and show that the topological invariance for the rotation number follows from this topologically flavored property.  Fortunately, area preserving and irrationally indifferent holomorphic maps enjoy this property.  In \cite{GAPELE96}, Gambaudo-Le Calvez-P\'ecou  show that a non linear rotation either verifies the linking property or a Birkhoff--P\'erez-Marco property associated to the existence of completely invariant non-trivial compact sets (the so called $\mathcal{P}$ condition). Further, by using the action on the prime ends they show that the Naishul's result holds under this alternative condition.  Hence, the panorama for the invariance of the rotation number is fairly complete in the surface setting. Coming back to the holomorphic case, it is easy to note that a rationally indifferent map do not verify the linking property. However, the  $\mathcal{P}$ condition holds and the Naishul's result follows. Nevertheless, the Leau-Fatou flower Theorem asserts that the local dynamics is combinatorially finite and quite simple. Hence, the invariance of the rotation number is direct. Recently, Le Roux \cite{LERO08} show in general that for any homeomorphic non linear rotation this more simple alternative occurs. Indeed, Le Roux result says that in the absence of the linking property,  the map is topologically conjugated to a rationally indifferent polynomial and a topological Leau-Fatou flower appears. Summarizing, in the surface setting, the Naishul's result holds either by the linking property or by the existence of a Leau-Fatou-Le Roux flower, and the prime ends machinery can be bypassed. \\

In higher dimensions, Gambaudo and P\'ecou \cite{GAPE95} treat the case of a differentiable local homeomorphism  admitting an invariant real codimension two torus. They also assume that the restricted dynamics on this torus is conjugated to an irrational translation. They show that the complementary tangent direction wraps around this invariant torus with a well defined asymptotic speed defining a rotation number. Provided the linking property holds, this rotation number is an invariant under topological conjugacies. \\

In this work we deal with a similar situation. We consider fibred holomorphic maps over irrational rotations. Given an indifferent invariant curve, a fibred rotation number is computed as the mean of the rotation angles at each fibre. We show that this fibred rotation number is an invariant under fibred topological conjugacies. There are some important remarks to do:  at this time,  we don't know if the linking property holds in this situation; on the other hand, there is not a well understood theory of {\it parabolic}  behavior for fibred maps and so, we don't have an alternative {\it \`a la Le Roux}. Even though there is no  satisfactory high dimensional prime ends theory, the skew product structure allow us to  use the one dimensional prime ends  technology in this case.   \\

\noindent 
$\mathbf{Acknowledgements. } $ I want to thank Jean-Marc Gambaudo for very useful conversations and comments. Also, I must acknowledge the sympathy, patience and guidance of Patrice Le Calvez in face of many wrong comments I made regarding the prime ends technique applied to this problem.   
\section{Fibred holomorphic maps}
Let $\alpha$ be an irrational angle in $\T$. A {\it fibred holomorphic map} is a continuous transformation
\begin{eqnarray*} 
F:\T\times \D &\longrightarrow& \T\times \C\\
 (\theta,z)&\longmapsto& \big( \theta+\alpha,f_{\theta}(z)\big)
 \end{eqnarray*}
and such  that the functions $f_{\theta}:\D\to \C$ are univalent for all $\theta\in \T$. We also assume that the zero section $\T\times \{0\}$ is an invariant curve, that is,  $f_{\theta}(0)=0$ for every $\theta\in \T$. Invariant curves  play the role of a center around which  the dynamics of $F$ is organized, generalizing thus the role of a fixed point for the local dynamics of an holomorphic germ (see \cite{PONC07}). 
We say that the invariant curve is {\it indifferent} if 
\begin{equation*}
\int_{\T}\log \big|\partial_zf_{\theta}(0)\big|d\theta=0.
\end{equation*} 
We recall that $F$ is injective and so the differential $\partial_zf_{\theta}$ is always non zero. As shown in \cite{PONC07}, a non indifferent invariant curve is either attracting or repelling in the sense that there exists a topological tube which is attracted (or repelled) to the curve in the future. In fact, this is an equivalent definition for being  indifferent. 
 Let's suppose  the topological degree of the application $\theta\mapsto \partial_zf_{\theta}(0)$ is zero. In other words, the  application  $\theta\mapsto \partial_{z}f_{\theta}(0)$
 is homotopic  in $\C\setminus \{0\}$ to a constant. Under this hypothesis we can define the logarithm of $\partial_z f_{\theta}(0)$.  We define a number which represents  the average rotation speed of the dynamics around the invariant curve, 
\begin{equation*}
\varrho_{T}(F)=\frac{1}{2\pi i }\int_{\T}\log \partial_zf_{\theta}(0)d\theta .
\end{equation*}
This is a real number and is called  the \emph{fibred rotation number}.  Note the $\log$ above is well defined $\mod 2\pi i$ and hence the number $\varrho_{T}(F)$ is well defined $\mod \Z$.  
In \cite{PONC07} the author also shows that in the analytic class of fibred holomorphic maps, and under an extra diophantine condition on the pair $(\alpha, \varrho_T(F))$, the map is conjugated to the pure linear map $(\theta, z)\mapsto (\theta+\alpha, e^{2\pi i \varrho_T(F)}z)$. 
\\

Let's study the behavior of this rotation number under conjugacies. 
A continuous map $f:U\subset(\C,0)\to \C$  defined on a neighborhood of the origin will be called a {\it local homeomorphism} if it is a orientation preserving homeomorphism  onto its image and left fixed the origin.  We will  consider continuous change of coordinates $H$  defined on a tubular neighborhood of the invariant curve in the form
\begin{equation*}
(\theta,z)\stackrel{H}{\longmapsto}\big(\theta,h_{\theta}(z)\big),
\end{equation*}
where  functions $h_{\theta}$ are local homeomorphisms. We call such a conjugacy a {\it fibred conjugacy}. When conjugating the map $F$ we will get a new fibred map $\tilde{F}=H^{-1}\circ F\circ H$ having the zero section  as an invariant curve. The topological characterization of an indifferent invariant curve implies this curve is also  indifferent.  Suppose  for a while that $h_{\theta}$ is differentiable at the origin and $\theta\mapsto \partial_z h_{\theta}(0)$ has zero degree. An easy computation implies that the invariant curve has zero degree and the fibred rotation number does not change, that is, $\varrho_{T}(F)=\varrho_T(\tilde{F})$. 
\\

As discussed in the Introduction, we are interested in to show that the fibred rotation number is a topological characterization of local dynamics, and not just a differentiable one. Hence, we are interested in the invariance of the fibred rotation number under fibred continuous conjugacies. The topological degree of the derivative at the origin need a topological counterpart: we require that $H$ is isotopic to the identity.  The remaining part of this work is devoted to prove the
\begin{theo}[Naishul's Theorem for fibred holomorphic maps]
Let $F$ and $G$ be two fibred holomorphic maps, admitting the zero section as an indifferent, zero degree invariant curve. Suppose $F$ and $G$ are topologically conjugated by a fibred conjugacy isotopic to the identity. Then the fibred rotation numbers  $\varrho_{T}(F)$, $\varrho_T(G)$ are equals.
\end{theo}
The structure of the paper and the proof of the main theorem is as follows: in Section 3 we define a rotation number for a finite fibres circle homeomorphism and we show that this number is a topological conjugacy invariant.  In Section 4 we define a rotation number for a finite fibres dynamics on compact sets by the prime ends action and the results on Section 3. Then we show that the fibred rotation number of the invariant curve may be computed as the limit of rotation numbers associated to finite fibres dynamics on compact sets. As the finite fibres rotation number is a topological conjugacy invariant, the limit rotation number is also invariant and the theorem holds.
 
\section{Finite fibres circle homeomorphisms}
Let $g_0, g_1, \dots , g_{n-1}$ be positive circle homeomorphisms. We define the finite fibres circle homeomorphism $G$ by
\begin{eqnarray*}
G:\Z_n\times \T&\longrightarrow& \Z_n\times \T\\
(i,x)&\longmapsto& \big(i+1,g_i(x)\big).
\end{eqnarray*}
The case $n=1$ reduces to the well known case of a single positive circle homeomorphism. Since Poincar\'e one knows that for this case, all orbits are described (at least combinatorially) by a unique real number, the so called {\it rotation number}.\\
 Let us denote by $\tilde g_i:\R\to \R$ a lift of $g_i$ to the real line, that is:
 \begin{itemize}
 \item[i) ] $\tilde g_i(x+1)=\tilde g_i(x)+1.$
 \item[ii) ] $\Pi\circ \tilde g_i=g_i\circ \Pi$, where $\Pi:\R\to \T$ stand for the natural projection.
 \end{itemize}  
 These lifts are continuous, increasing and define a finite fibres dynamics 
\begin{eqnarray*}
\tilde G:\Z_n\times \R&\longrightarrow& \Z_n\times \R\\
(i,x)&\longmapsto& \big(i+1,\tilde g_i(x)\big).
\end{eqnarray*}
In order to take account of the rotation speed of orbits we define the {\it $m^{th}$-step} of a point $(i,x)\in \Z_n\times \R$ by the function
\[
\Psi^{(m)}(i,x)=\Pi_\R\tilde G^{m}(i,x)-x.
\]
where $\Pi_\R:\Z_n\times \R\to \R$ stands for the second coordinate projection.
The following is an easy consequence of the definition of $\tilde G$ (the reader can refer to  \cite{FEGLPIST02}):
\\

\begin{lema}\label{sec2p1}\noindent
\begin{itemize}
\item[i) ] $\psi^{(m)}$ is periodic on $x$, that is, for every $k\in \Z$, $x\in \R$ and $i\in \Z_n$ we have $$\Psi^{(m)}(i,x+k)=\Psi^{(m)}(i,x).$$
\item[ii) ] For every $x,x'\in \R$, $i\in \Z_n$ we have 
\[
|\Psi^{(m)}(i,x)-\Psi^{(m)}(i,x')|\leq 1\quad_{\blacksquare}
\]
\end{itemize}
\end{lema}
Previous Lemma will allow us to show that the sequence $\{\Psi^{(m)}(i,x)\}_{m\in \N}$ verifies a sub-additive property.
\begin{lema}
 There exist a constant $C>0$ such that for every $x\in \R$, $i\in \Z_n$, $m,r\in \N$ one has
\[
\Psi^{(m+r)}(i,x)\leq \Psi^{(m)}(i,x)+\Psi^{(r)}(i,x)+C.
\]
\end{lema}
{\it Proof. } 
\begin{eqnarray*}
\Psi^{(m+r)}(i,x)&=&\tilde G^{m}\big(i+r,\tilde G^{r}(i,x)\big)-\tilde G^r(i,x)+\tilde G^r(i,x)-x\\
&=&\Psi^{(m)}\big(i+r, \tilde G^r(i,x)\big)+\Psi^{(r)}(i,x)\\
&\leq& \Psi^{(m)}(i+r, x)+\Psi^{(m)}(i,x)+1.
\end{eqnarray*}
For the last inequality we have used point $(ii)$ of Lemma \ref{sec2p1}. Note if $r=0(\Z_n)$ the we have the required inequality. In other case we proceed as follows: let $s\in \N$ such that $r+s=0(\Z_n)$. We have
\begin{eqnarray*}
\Psi^{(m)}(i+r,x)&\leq&s+\Psi^{(m+s)}(i+r,x)\\
&=&s+\Psi^{(m)}\big(i+r+s, \tilde G^s(i+r,x)\big)+\Psi^{(s)}(i+r,x)\\
&\leq& s+\Psi^{(m)}(i, x)+\Psi^{(s)}(i+r,x)+1.
\end{eqnarray*}
Hence, we can pick $C=\sup_{(s,x)\in \Z_n\times \R}s+\Psi^{(m)(n-s, x)}+2\quad_{\blacksquare}$
\begin{prop}
For every $i\in \Z_n$,  $x\in \R$ the limit
\[
\rho_{ff}(\tilde G, i, x)=\lim_{m\to \infty}\frac{\Psi^{(m)}(i,x)}{m}
\]
exists and belongs to $\R$. Moreover, this limit is independent of the choice $(i,x)$.
\end{prop}
{\it Proof. } The existence of the limit is guaranted by the sub-aditive Lemma (see for instance \cite{KAHA95}). Further, we know that $\Psi^{(m)}(i,x)+m>0$ which implies $\rho(\tilde G, i,x)>-\infty$.  The independence on $x$ is a consequence of the point $(ii)$ of Lemma \ref{sec2p1}.  Finally, for $i,j\in \Z_n$ and $m\in \N$ we have
\begin{eqnarray*}
\rho_{ff}(\tilde G, j, x)&=&\lim_{m\to \infty}\frac{\Psi^{(m)}(j,x)}{m}\\ &=&\lim_{m\to \infty}\frac{\Psi^{(m-(i-j))}\big(i,\tilde G^{i-j}(x)\big)}{m}+\lim_{m\to \infty}\frac{\Psi^{(i-j)}(j,x)}{m}\\
&=&\rho_{ff}(\tilde G,i,x)\quad_{\blacksquare}
\end{eqnarray*}
We define the {\it finite fibres rotation number of $\tilde G$} by a (any) rotation number $\rho_{ff}(\tilde G)=\rho_{ff}(\tilde G, i, x )$. The following properties are easy to show:
\begin{lema}
\noindent
\begin{itemize}
\item[i) ] Let $\hat G$ be another lift  of  $G$. Then there exists an integer $p$ such that $\rho_{ff}(\tilde G)=\rho_{ff}(\hat G)+p$.
\item[ii) ] For every $m\in \Z$ one has $\rho_{ff}(\tilde G^{m})=m\rho_{ff}(\tilde G)\quad_{\blacksquare}$
\end{itemize}
\end{lema}
First point above allow us to define the rotation number for $G$ as 
\[
\rho_{ff}(G)=\rho_{ff}(\tilde G)  \mod (\Z).
\]
Moreover, we can always pick a lift $\tilde G$ so that $\rho_{ff}(\tilde G)\in [0,1)$.
\begin{rema}
The map $G^{n}\big|_{\{0\}\times \T}:\T\to \T$ is a positive circle homeomorphism. Hence, the asymptotic turning speed for his orbits is controlled by its  rotation number $\rho(G^{n}\big|_{\{0\}\times \T})$. This number coincides with the rotation number $\rho_{ff}(G^n)$ and hence we get
\begin{equation}
\rho(G^n\big|_{\{0\}\times \T})=n\rho_{ff}(G)\mod (\Z).\label{sec2dos}
\end{equation}
\end{rema}
\paragraph{The finite fibres rotation number is a fibred conjugacy invariant. }
We say that a continuous map $H:\Z_n\times \T\to \Z_n\times \T$ is a  fibred conjugacy if 
\[
H(i,x)=\big(i, h_i(x)\big)
\]
and $h_i:\T\to \T$ is a positive circle homeomorphism for every $i\in \Z_n$. Given a finite fibres circle homeomorphism $G:\Z_n\times \T\to \Z_n\times \T$, the conjugacy of $G$ by $H$ results in a new finite fibres circle homeomorphism $F=H^{-1}\circ G\circ H$. The rotation number is a topological invariant under this kind of conjugacy:
\begin{prop}
With the previous notations, we have $\rho_{ff}(F)=\rho_{ff}(G)$.
\end{prop} 
{\it Proof. } Let $ \tilde G$ be a lift of $G$ and $\tilde H$ a lift of $H$ so that $\Pi_{\R}\tilde H(i,0)\in [0,1)$ for every $i\in \Z_n$.  Then $\tilde F=\tilde H^{-1}\circ \tilde G\circ \tilde H$ is a lift for $F$. By considering the $n^{th}$ iterate of both, $\tilde G$ and $\tilde F$, we get two topologically conjugated positive circle homeomorphism on the first fibre. By using (\ref{sec2dos}) and the topological invariance of the rotation number for single homeomorphisms we have
\begin{eqnarray*}
\rho(\tilde G^n\big|_{\{0\}\times \R})&=&\rho(\tilde F^n\big|_{\{0\}\times \R})\\
n\rho_{ff}(\tilde G)&=&n\rho_{ff}(\tilde F)\\
\rho_{ff}( G)&=&\rho_{ff}( F) \quad_{\blacksquare}
\end{eqnarray*}
\section{The rotation number for Birkhoff--P\'erez-Marco's invariant compact sets}\label{tres}
In this section we will assign a rotation number to certain classes of completely invariant connected compact sets. We start by a classical construction due to Cartwright and Littlewood \cite{CALI51} that assigns a rotation number to a completely invariant full compact subset of $\C$ under the action of a  local homeomorphism.  
\\

Let $K\subset \overline{\C}$ be a nontrivial compact connected full set (has more than one point and it complement is simply connected). We assume $0\in K$.  By the Riemann mapping theorem, there exists a biholomorphic map $\phi: \overline{\C}\setminus \D \to \overline{\C}\setminus K $, normalized so that $\phi$ is tangent to the identity at $\infty$. In general this map do not extend to the boundary. Cartwright and Littlewood introduce the topological compact space of {\it prime-ends} $\mathcal{P}(K)$ which is homeomorphic to the circle $\T$ and  such that $\phi$ extends continuously  to the boundary $\partial \D$ as a map onto the compactification of   $ \overline{\C}\setminus K$ by $\mathcal{P}(K)$. Nowadays this is known as the {\it prime-ends compactification}. 

 Let $F$ be a local homeomorphism defined on a neighborhood of $K$ and fixing the origin. Let's assume that $K$ is completely invariant by $F$, that is $F(K)=F^{-1}(K)=K$. The map $F$ extends to the prime ends compactification as a homeomorphism, and its action on $\mathcal{P}(K)$ is a positive circle homeomorphism. We denote by $\rho(F,K)$ the rotation number associated to this circle homeomorphism.  This rotation number is a topological conjugacy invariant.
 
 In order to compute this rotation number, we need to recall that this number takes account of the orbits rotation speed of any point in $\T$ (or $\mathcal{P}(K)$). This speed can be computed using the advancing of a point in $\T\cong \mathcal{P}(K)$ under the action of a lift of $F$ to the real line. 
 
 More precisely, let $F_{\T}:\R\to \R$ be a lift of the action of $F$ on the prime ends  $\T\cong \mathcal{P}(K)$ and $x\in \T$. The classical theory of the rotation number gives
 \begin{prop}
 \noindent
 \begin{itemize}
 \item[i) ] The limit \[
 \lim_{n\to \infty}\frac{F_{\T}^n(x)-x}{n}
 \]
 exists and, up to an integer, equals $\rho(F,K)$. 
 \item[ii) ] The ordering of the iterates $\{F_{\T}^n(x)\}_{n\in \N}$ with respect to the $\Z$ lattice determines uniquely the rotation number$\quad_{\blacksquare}$
  \end{itemize}
 \end{prop}
 Hence, in order to compute the rotation number $\rho(F,K)$ just by using the map $F$, we need to identify a single prime end and keep track of the ordering of its orbit on $\mathcal{P}(K)$. In general a single prime end is not represented  by a single point on $K$. The case is simpler when a point $\tilde z\in K$ is accessible, that is, there exists a path $\gamma\in \overline{C}\setminus K$ such that $\tilde z$ is an extremal point of $\gamma$.  In this case we can assume that $\tilde z$ is a prime end of $K$.   The map $F$ acts on $\C\setminus \{0\}$, and we can lift $F$ to the universal covering $\R\times \R^{+}$ (the upper half plane). We continue to denote by $F$ this lifting and by $K$ the lifting of $K\setminus{0}$. The previous discussion allows us to state 
 \begin{prop}
  Let $\Pi_{\R}:\R\times \R^{+}\to \R$ be the canonical projection. Let $\tilde z\in K$ be an accessible point. Then the limit
  \[
 \lim_{n\to \infty}\frac{\Pi_{\R}\left(F^n(\tilde z)-\tilde z\right)}{n}
 \]
exists and, up to an integer, equals $\rho(F, K)\quad_{\blacksquare}$
 \end{prop}
The type of local homeomorphisms $F$ we have in mind are indifferent holomorphic germs in the plane. Let $F(z)=e^{2\pi i \beta}z+r_2z^2+r_3z^3+\dots$ be a holomorphic map. The number  $\beta\in \T$ is called the {\it rotation number} of $F$. In \cite{PERE97} P\'erez-Marco shows that there exists a completely invariant, nontrivial connected compact  full set  $K$ containing the fixed point $z=0$. Moreover we have
\begin{prop}[P\'erez-Marco \cite{PERE97}]
The rotation number $\rho(F,K)$ equals $\beta$. Furthermore, $\beta$ is the rotation number of any completely invariant nontrivial compact connected full set   containing the fixed point 
$z=0\quad_{\blacksquare}$ 
\end{prop}
As a direct consequence of this Proposition and the topological invariance of the rotation number on the prime ends, P\'erez-Marco  gives a new proof to the  corresponding Naishul's Theorem. That is, the rotation number  is a topological conjugacy invariant for indifferent holomorphic germs in the plane. 
\subsection{Rotation number for a finite fibres Birkhoff--P\'erez-Marco's invariant}\label{sec31}
Let $n\in \N$. We say that a $n-tuple$ of local homeomorphism  $(f_0,f_1,\dots,f_{n-1})$  is a \emph{good chain} for the $n-tuple$ of sets $(U_0,U_1,\dots,U_{n-1})$ if for every $i\in \Z_n$ we have
 \begin{itemize}
 \item[$i)$ ] $U_i$ is an open Jordan domain (the interior of a Jordan curve) containing the origin in $\C$.
 \item[$ii)$ ] $f_i$ is defined and injective in a neighborhood of  $\overline{U_i}$ and $f_i^{-1}$ is a local homeomorphism defined and injective in a neighborhood of  $\overline{U_{i+1}}$.
 \end{itemize} 
 We say that a good chain as above has a {\it finite fibres Birkhoff--P\'erez-Marco compact set} if there exists compact connected full sets $K_i\subset \overline{U_i}$,  containing the origin but not reduced to it, and such that 
      \begin{equation*}
      f_i(K_i)=K_{i+1}\quad , \quad f_{i}^{-1}(K_{i+1})=K_i 
      \end{equation*} for all  $i$ in $\Z_n$.    The good chain and the Birkhoff-P\'erez-Marcos's invariant compact set give raise to  an action
      \begin{eqnarray*}
      G:\bigcup \mathcal{P}(K_i)&\longrightarrow& \bigcup\mathcal{P}(K_i)\\
      (i,x)&\longmapsto&\big(i+1, g_i(x)\big),
      \end{eqnarray*}
where $g_i:\T\to \T$ is the positive circle homeomorphism given by the action of $f_i$ from  $\mathcal P(K_i)$ onto  $\mathcal{P}(K_{i+1})$.  Note that the normalization on the choice of the Riemmann mapping is used here in order to identify a common origin for each circle. In this way, $G$ is a finite fibre circle homeomorphism and the rotation number $\rho_{ff}(G)$ can be computed by means of the equality $n\rho_{ff}(G)=\rho(G^n\big|_{\mathcal{P}(K_0)})$.      By putting 
\[
F=f_{n-1}\circ f_{n-2}\circ \dots \circ f_0
\]
we get that $G^n\big|_{\mathcal{P}(K_0)}$ corresponds to the action of $F$ on the prime-ends of the completely invariant compact set $K_0$. 
\\

When the functions $f_i$ are holomorphic and the fixed point $z=0$ is indifferent, that is, $|\partial_zF(0)|=1$, 
a Theorem by P\'erez-Marco (see \cite{PERE97, PONC07}) guaranties the existence of a finite fibres  compact set.  Under this holomorphic hypothesis we can pick a lift $\tilde F$ such that
 \[
\rho(\tilde F,K_0)=\frac{1}{2\pi i}\sum_{i=0}^{n-1}\log f_i'(0).
\]
This lift induces a lift $\tilde G$ and the following holds:
\begin{lema}\label{lema10}
The rotation number for the finite fibres circle homeomorphism $G$ induced by the good chain of holomorphic local homeomorphisms  $(f_0,f_1,\dots,f_{n-1})$  is
\[
\rho_{ff}(\tilde G)=\frac{1}{2\pi i n}\sum_{i=0}^{n-1}\log f_i'(0)\quad_{\blacksquare}
\]
\end{lema}
\subsection{Rotation number for a fibred P\'erez-Marco's invariant}
Consider the fibred holomorphic map
\begin{eqnarray*}
F:\T\times \D&\longrightarrow&\T\times \C\\
(\theta, z)&\longmapsto& \big(\theta+\alpha, \rho_1(\theta)z+\rho_2(\theta)z^2+\dots\big).
\end{eqnarray*} 
 Assume that the invariant curve $\T\times \{0\}$ has zero degree and is indifferent. Under this hypothesis we have
\begin{theo}[Fibred     P\'erez-Marco's invariant \cite{PONC07}]
There exists $r_0>0$ such that  for every $0<r\leq r_0$ the fibred holomorphic maps  $F:\T\times\D_r\to \T\times \C$ and $F^{-1}:\T\times\D_r\to \T\times \C$ are  local diffeomorphisms on each fibre. Moreover, there exists a compact connected compact set $K_r\subset \T\times \overline{\D_r}$ such that
\begin{itemize}
\item[i) ] $K_r$ contains the invariant curve $\T\times \{0\}$.
\item[ii) ] The fibres $K_{r,\theta}\subset \overline{\D_r}$ are  connected full compact sets.
\item[iii) ] $F(K_r)=F^{-1}(K_r)=K_r$.
\item[iv) ] There exists $\tilde \theta\in \T$ such that $K_{r,\tilde \theta}\cap \partial \D_{r}\neq \emptyset$.
\end{itemize}
\end{theo} 
This completely invariant compact set can be obtained as follows: fix a sequence $\{\frac{p_n}{q_n}\}_{\N}$ of rational numbers converging to $\alpha$.  Construct the finite approximations of $F$ as 
\begin{eqnarray*}
F_{(n)}:\T\times \D_r&\longrightarrow &\T\times \C\\
(\theta, z)&\longmapsto& \left(\theta+\frac{p_n}{q_n}, \rho_1^{(n)}(\theta)z+\rho_2(\theta)z^2+\dots\right)
\end{eqnarray*}
and such that $\rho_1^{(n)}\to \rho_1$ uniformly on $\T$. In this way, the maps $F_{(n)}\to F$ uniformly. Provided that $\sum_{0}^{q_n-1}\log |\rho_1^{(n)}(jp_n/q_n)|=0$, each dynamics $F_{(n)}$ has a compact set $K_{r}^{n}$ with the desired topological and invariant properties (a finite fibres invariant compact set). By taking a Haussdorf converging sub-sequence $K_r^{n}$, we get the desired invariant compact set $K_r$ (a filling process on each fibre may be also necessary). For details, see Lemma 4.20 in \cite{PONC07}. 
\\

This construction suggests a way for associating a rotation number to $K_r$. We say that $\beta \in \R$ is a {\it fibred rotation number} for $K_r$ if 
\begin{itemize}
\item[i) ] There exists a sequence of continuous fibred dynamics
\begin{eqnarray*}
F_{(n)}:\T\times \D_r&\longrightarrow &\T\times \C\\
(\theta, z)&\longmapsto& \left(\theta+\frac{p_n}{q_n},f_{n,\theta}(z)\right)
\end{eqnarray*}
where $f_{n,\theta}:\D_r\to \C$ is a local homeomorphism. We also require that the zero section $\T\times \{0\}$  is invariant by every $F_{(n)}$.
\item[ii) ] There exists compact sets $K^{n}\subset \T\times \D_r$ containing the zero section, such that
\begin{itemize}
\item[j) ] There exists $\tilde \theta\in \T$ such that the fibres $K^n_{ \theta}=\{0\}$ for every $\theta\in \T\setminus\{\tilde \theta+j \frac{p_n}{q_n}\}_{j\in \Z_n}$ and $K^n_{\tilde \theta}\neq \{0\}$. Moreover, these non-trivial fibres are compact, connected and full sets.
\item[jj) ] $F_{(n)}(K^{n})=F_{(n)}^{-1}(K^n)=K^n$.
\end{itemize}
Thus, each $K^n$ is a finite fibres compact set, completely invariant  by $F_{(n)}$ when restricted to the fibres $\{\tilde \theta+j \frac{p_n}{q_n}\}_{j\in \Z_n}$. A rotation number $\rho_{ff}(F_{(n)}, K^n)$ is associated for the action of $F_{(n)}$ on the prime ends as done in Section \ref{sec31}.   
\item[iii) ] $K$ is the Haussdorf limit of $K^n$, unless a fulling process on each fibre. 
\item[iv) ] $\beta=\lim_{n\to \infty}\rho_{ff}(F_{(n)}, K^n)$.
\end{itemize}
Note that in this definition we require that the sequences $K^n$ and $\rho_{ff}(F_{(n)}, K^n)$ are both convergents.
\begin{exam}
The indifferent zero degree  curve $\T\times\{0\}$ is invariant by $F$. As defined at the Introduction, this curve has a  fibred rotation number 
\[
\varrho_T(F)=\frac{1}{2\pi i}\int_{\T}\log \rho_1(\theta)d\theta.
\]
This number is a fibred rotation number for $K$. Indeed, the above construction for $K$ gives us the functions $F_{(n)}$ and the finite fibres compact invariants $K^n$. By Lemma \ref{lema10} the rotation number verifies 
\[
\rho_{ff}(\tilde F_{(n)}, K^n)=\frac{1}{2\pi i q_n}\sum_{j=0}^{q_n-1}\log \rho_1^{(n)}\left(\tilde \theta+j\frac{p_n}{q_n}\right).
\] 
This sum is a Riemman sum converging to the integral defining $\varrho_T(F)\quad_{\blacksquare}$
\end{exam}
\begin{lema}
The set of fibred rotation numbers for $K$ has at most one element.
\end{lema}
\noindent{\it Proof. } Let $( F_{(n)}, K^n)$ and $(G_{(n)}, L^n)$ two sequences of fibred dynamics and finite fibres compact set defining two rotation numbers $\beta_F, \beta_G$ for $K$. Let $\eps>0$ and two integers $r, s$ such that 
\[
\beta_F<\frac{r}{s}<\beta_F+\frac{\eps}{3}.
\]
Since $F_{(n)}$ and $G_{(n)}$ converge uniformly to $F$ we can pick $n^*\in \N$ and $\delta>0$ such that
\begin{equation*}
d(F_{(n^*)}^j(x,\theta), G_{(n^*)}^j(y,\theta))\leq \frac{\eps}{3}
\end{equation*}
for  every pair $x, y$ in $\C$ verifying  $d(x,y)<\delta$, every $\theta\in \T$ and $0\leq j\leq s$. The sequences $K^n, L^n$ Haussdorf converge to $K$, and $\rho_{ff}(F_{(n^*)}, K^{n^*})\to \beta_F$, hence we can also assume  that $$d_H(K^{n^*}, L^{n^*})<\delta\quad \textrm{and}\quad |\rho(F_{(n^*)} ,K^{n^*})-\beta_F|<\frac{\eps}{3}.$$  The accessible points are dense on the boundary of $K^{n^*}$ and $L^{n^*}$. Hence, we can pick accessible points $(\tilde \theta, \tilde x)\in K^{n^*}$ and $(\tilde \theta, \tilde y)\in L^{n^*}$ such that $d(\tilde x, \tilde y)<\delta$.  By computing the advancing of this accessible points in a suitable lifting (as done in Section \ref{tres}), we deduce 
\[
\beta_G<\beta_F+\eps.
\]
In an analogous way we deduce $\beta_F<\beta_G+\eps$ and the equality holds $\quad_{\blacksquare}$\\

As seen in the Example above, the number $\varrho_T(F)$ is actually a fibred rotation number for $K$ and hence is the only one.  By the topological nature of the construction and the unicity, this fibred rotation number is a topological invariant by fibred conjugacies, as shown by the following Proposition, which implies directly the Theorem 2:
\begin{prop}
Let $F, G$ be fibred holomorphic maps over the irrational rotation $\theta\mapsto \theta+\alpha$. Assume that the zero section is an invariant, indifferent and zero degree curve for both dynamics. Let $H(\theta, z)=(\theta, h_{\theta}(z))$  be a  isotopic to the identity fibred  conjugacy, such that each fibre map $h_{\theta}$ is a local homeomorphism fixing the origin. Let $K_F$ be a fibred P\'erez Marco's compact set for $F$ and $K_G=H^{-1}(K_F)$ the corresponding P\'erez-Marco's compact set for $G$.  Then the fibred rotation numbers are equals, that is, $\varrho_T(F)=\varrho_T(G)$.  
\end{prop}
\noindent
{\it Proof. }  Let $F_{(n)}$ and $K^n$ the sequences of fibred maps and finite fibres compact sets used in the construction of $K_F$. The sequences
\[
G_{(n)}=H^{-1}\circ F_{(n)}\circ H\quad\textrm{and}\quad L^{n}=H^{-1}(K^n)
\]
are corresponding sequences converging uniformly to $G$ and defining the invariant compact set $K_G=H^{-1}(K_F)$. 
As $H$ is isotopic to the identity, we can pick a lift $\tilde H$ such that $\Pi_{\R}\tilde h_{\theta}(0)\in [0,1)$ for every $\theta \in \T$ (here $0\in \T$ means for the origin in each prime ends circle, given by the real direction in $\C$). 
By definition, the finite fibres dynamics $\tilde F_{(n)}$ and $\tilde G_{(n)}$ are topologically conjugated, and thus the rotation numbers of the action over the prime ends of $K^n$ and $L^n$ are the same, that is, $\rho_{ff}(F_{(n)}, K^n)=\rho_{ff}(G_{(n)}, L^n)$. This sequence of real numbers converges in one hand to the unique fibred rotation number of $K_F$, that is, , $\varrho_T(F)$, and in the other, to the unique fibred rotation number for $K_G$, namely  $\varrho_T(G)$,  and the equality holds $\quad_{\blacksquare}$
\bibliographystyle{plain}
  \bibliography{naishul.bib}

\def\cprime{$'$}
\begin{thebibliography}{1}

\bibitem{CALI51}
M.~L. Cartwright and J.~E. Littlewood.
\newblock Some fixed point theorems.
\newblock {\em Ann. of Math. (2)}, 54:1--37, 1951.
\newblock With appendix by H. D. Ursell.

\bibitem{GAPELE96}
Jean-Marc Gambaudo, Patrice Le~Calvez, and {\'E}lisabeth P{\'e}cou.
\newblock Une g\'en\'eralisation d'un th\'eor\`eme de {N}aishul.
\newblock {\em C. R. Acad. Sci. Paris S\'er. I Math.}, 323(4):397--402, 1996.

\bibitem{GAPE95}
Jean-Marc Gambaudo and Elisabeth P{\'e}cou.
\newblock A topological invariant for volume preserving diffeomorphisms.
\newblock {\em Ergodic Theory Dynam. Systems}, 15(3):535--541, 1995.

\bibitem{KAHA95}
Anatole Katok and Boris Hasselblatt.
\newblock {\em Introduction to the modern theory of dynamical systems},
  volume~54 of {\em Encyclopedia of Mathematics and its Applications}.
\newblock Cambridge University Press, Cambridge, 1995.
\newblock With a supplementary chapter by Katok and Leonardo Mendoza.

\bibitem{LERO08}
Fr{\'e}d{\'e}ric Le~Roux.
\newblock A topological characterization of holomorphic parabolic germs in the
  plane.
\newblock {\em Fund. Math.}, 198(1):77--94, 2008.

\bibitem{NAIS82}
V.~A. Na{\u\i}shul{\cprime}.
\newblock Topological invariants of analytic and area-preserving mappings and
  their application to analytic differential equations in {${\bf C}\sp{2}$} and
  {${\bf C}P\sp{2}$}.
\newblock {\em Trudy Moskov. Mat. Obshch.}, 44:235--245, 1982.

\bibitem{PERE97}
Ricardo P{\'e}rez-Marco.
\newblock Fixed points and circle maps.
\newblock {\em Acta Math.}, 179(2):243--294, 1997.

\bibitem{PONC07}
Mario Ponce.
\newblock Local dynamics for fibred holomorphic transformations.
\newblock {\em Nonlinearity}, 20(12):2939--2955, 2007.

\bibitem{FEGLPIST02}
J.~Stark, U.~Feudel, P.~A. Glendinning, and A.~Pikovsky.
\newblock Rotation numbers for quasi-periodically forced monotone circle maps.
\newblock {\em Dyn. Syst.}, 17(1):1--28, 2002.

\end{thebibliography}

\end{document}